\documentclass{article}

\usepackage{latexsym}
\usepackage{a4wide}
\usepackage{amscd}
\usepackage{graphics}
\usepackage{amsmath}
\usepackage{amssymb}
\usepackage{mathrsfs}
\input xy
\xyoption{all}

\newcommand{\N}{\mathbb{N}}                     
\newcommand{\Z}{\mathbb{Z}}                     
\newcommand{\R}{\mathbb{R}}                     
\newcommand{\set}[2]{\left\{{#1}\mid{#2}\right\}}       
\newcommand{\qed}{\hfill $\Box$ \bigskip}       
\newcommand{\dist}{\mathrm{dist\,}}             

\newtheorem{thm}{\sc Theorem}[section]      
\newtheorem{cor}[thm]{\sc Corollary}        
\newtheorem{lem}[thm]{\sc Lemma}            
\newtheorem{prop}[thm]{\sc Proposition}     
\newtheorem{defn}[thm]{\sc Definition}      
\newtheorem{rem}[thm]{\sc Remark}	    

\title{High action orbits for Tonelli Lagrangians and superlinear
  Hamiltonians on compact configuration spaces}
\author{Alberto Abbondandolo\thanks{
    Dip. di Mat.,
    Universit\`a di Pisa,
    Largo Bruno Pontecorvo 5, 56127 Pisa, Italy.
        \textsf{e-mail: abbondandolo@dm.unipi.it}},\ \ \
    Alessio Figalli\thanks{
    Scuola Normale Superiore,
    Piazza  dei Cavalieri 7, 56126 Pisa, Italy.
        \textsf{e-mail: a.figalli@sns.it}}}
\date{September 22, 2006}

\begin{document}

\maketitle

\begin{abstract}
Multiplicity results for solutions of various boundary value problems
are known for dynamical systems on compact configuration manifolds, 
given by Lagrangians or Hamiltonians which have quadratic growth in 
the velocities or in the momenta. 
Such results are based on the richness of the topology of the space of
curves satisfying the given boundary conditions.
In this note we show how these results can be extended to the classical
setting of Tonelli Lagrangians (Lagrangians which are $C^2$-convex and
superlinear in the velocities), or to Hamiltonians which are
superlinear in the momenta and have a coercive action integrand. 
\end{abstract}  

\section*{Introduction}

Let $M$ be a compact manifold, the configuration space of a -
possibly non-autonomous - Lagrangian system, given by
a regular function $L$ - the Lagrangian - on $\R \times TM$. Here
$TM$ is the tangent bundle of $M$, and its elements are denoted by
$(q,v)$, where $q$ is a point in $M$ and $v$ is a tangent vector
at $q$. A classical assumption is that $L$ should be a {\em
Tonelli Lagrangian}, meaning that $L$ is fiberwise $C^2$-strictly
convex (that is $\partial_{vv}L>0$), and has superlinear growth on
each fiber.

On one hand Tonelli assumptions imply that the Legendre transform
is a diffeomorphism between the tangent bundle $TM$ and the
cotangent bundle $T^*M$, so that the Hamiltonian vector field on
$T^*M$ associated to the Fenchel transform of $L$ induces a vector
field on $TM$. Together with the assumption that such a vector
field is complete\footnote{This is always true if $L$ is
autonomous, due to the conservation of energy, to the coercivity
of the Hamiltonian, and to the compactness of $M$. In general,
completeness follows from a growth condition on $\partial_t L$.
See section \ref{bflt}.}, this produces a dynamical system on
$TM$.

On the other hand Tonelli assumptions imply existence results for
minimal orbits, such as the existence of an orbit having minimal
action connecting two given points on $M$, or the existence of a
periodic orbit of minimal action in every conjugacy class of the
fundamental group of $M$ (if $L$ is assumed to be periodic in
time). See for instance \cite{man91}.

However, more sophisticated existence results for orbits with
large action and large Morse index are known only for a smaller
class of Lagrangians, namely Lagrangians which behave
quadratically in $v$ for $|v|$ large. Their proofs are based on
topological methods in critical point theory. The fact that $L$
grows {\em at most} quadratically allows to find a nice functional
setting where the action functional is regular; the fact that $L$
grows {\em at least} quadratically implies the Palais-Smale
condition, the standard compactness property which makes minimax
arguments work. For instance, in the case of periodic Lagrangians
in this class, Benci \cite{ben86} has proved the existence of
infinitely many contractible periodic orbits, provided that $M$
has finite fundamental group.

The aim of this note is to show how the existence and multiplicity
results holding for Lagrangians which behave quadratically in $v$
for $|v|$ large can actually be proved for any Tonelli Lagrangian
which induces a complete vector field. We deal with many boundary
conditions at once by fixing a submanifold $Q$ of $M\times M$ and
by considering orbits $\gamma:[0,1]\rightarrow M$ satisfying the
non-local boundary condition
\begin{equation}
\label{uno} \begin{split}(\gamma(0),\gamma(1)) & \in  Q, \\
D_v L(0,\gamma(0),\gamma'(0))[\xi_0] & = D_v
L(1,\gamma(1),\gamma'(1))[\xi_1],\end{split}
\end{equation}
for every vector $(\xi_0,\xi_1)$ in $T(M\times M)$ which is
tangent to $Q$ at $(\gamma(0),\gamma(1))$. When $Q$ is the
diagonal submanifold, (\ref{uno}) yields to periodic orbits; when
$Q$ is the singleton $\{(q_0,q_1)\}$, (\ref{uno}) yields to orbits
joining $q_0$ to $q_1$. Other choices of $Q$ allow to deal with
orbits which are suitably normal to two submanifolds at their
end-points, or satisfy Neumann boundary conditions.

Our main result here is that the number of orbits satisfying
(\ref{uno})of a complete Tonelli Lagrangian system is greater than
the cuplength of the space $C_Q([0,1],M)$ consisting of all
continuous paths $\gamma:[0,1] \rightarrow M$ such that
$(\gamma(0),\gamma(1))$ belongs to $Q$ (the cuplength of a
topological space $X$ is the maximum length of a non-vanishing cup
product of elements of degree at least one in the homology ring of
$X$). Moreover, if $C_Q([0,1],M)$ has infinitely many
non-vanishing Betti numbers\footnote{It may happen that
$C_Q([0,1],M)$ has infinitely many non-vanishing Betti numbers but
its cuplength is finite, as in the case of the free loop space of
the sphere, see \cite{sul75}.}, there
is a sequence of such orbits with diverging action. See section
\ref{lagstat} for precise statements and for their consequences.
By taking a manifold $M$ with finite fundamental group and by
choosing $Q$ to be the diagonal in $M\times M$, we obtain the
generalization of Benci's theorem to Tonelli Lagrangians.

The idea of the proof is to modify the Tonelli Lagrangian for
$|v|$ large, making it quadratic in $v$ there. A suitable a priori
estimate on Lagrangian orbits with bounded action then allows to
prove that orbits of the modified Lagrangian found by a minimax
argument must lie in the region where the Lagrangian is not
modified. The result about the existence of infinitely many orbits
uses also Morse index estimates.

The same idea works for first order Hamiltonian systems on the
cotangent bundle $T^*M$ of the compact manifold $M$. The elements
of $T^*M$ are denoted by $(q,p)$, where $q$ is a point in $M$ and
$p$ is a cotangent vector at $q$. Here a dynamical system is
defined by a Hamiltonian $H:\R \times T^*M \rightarrow \R$. If we
do not assume that $H$ is convex on the fibers (or that
$\partial_{pp}H$ is everywhere invertible), such a Hamiltonian
system does not admit a Lagrangian formulation. In this setting,
Cieliebak \cite{cie94} has proved the analogue of Benci's result
for Hamiltonians which behave quadratically in $p$ for $|p|$ large
and whose action integrand $DH[Y]-H$ grows at least quadratically
in $p$ (here $Y$ is the Liouville vector field, expressed in local
coordinates by $\sum_i p_i \partial_{p_i}$). His
proof is based on the study of the Floer equation, a zero-order
perturbation of the equation for $J$-holomorphic curves on $T^*M$.
The quadratic behavior of $H$ and the growth condition on the
action integrand are crucial in the proof of the $L^{\infty}$
estimates for the solutions of the Floer equation having bounded
action.

Here we consider orbits $x:[0,1]\rightarrow T^*M$ satisfying the
general non-local boundary condition, which in the Hamiltonian
setting is just
\begin{equation}\label{tre} (x(0),-x(1))\in N^* Q,
\end{equation}
where $N^*Q$ is the conormal bundle of $Q$ in the cotangent bundle
of $M\times M$. We prove that if the Hamiltonian $H$ has
superlinear growth in $p$, its action integrand $DH[Y]-H$ is
coercive, and the Hamiltonian vector field is complete, then the
number of orbits satisfying (\ref{tre}) is greater than the
$\Z_2$-cuplength of $C_Q([0,1],M)$. Moreover, if infinitely many
$\Z_2$-homology groups of $C_Q([0,1],M)$ are non-trivial, then
there is a sequence of orbits satisfying (\ref{tre}) with
diverging action. When $Q$ is the diagonal of $M\times M$, we
obtain the generalization of Cieliebak's theorem to a larger class
of Hamiltonians. Notice also that if $L$ is a Tonelli Lagrangian
then its Fenchel transform $H$ is in the above class, so the
Hamiltonian setting considered here includes the Lagrangian one as
a particular case. Since the whole argument in the Lagrangian
setting is considerably simpler, we prove the Lagrangian
statements independently.

We wish to stress the fact that the quadratic modification
argument and the a priori estimate presented here seem to be quite
a general tool, and they should allow to generalize to Tonelli
Lagrangians or to the superlinear Hamiltonians in the above class
many other existence results, such as for instance Long's theorem
on the existence of infinitely many contractible orbits of integer
period for a time-periodic Lagrangian system on the torus, see
\cite{lon00b}.

\section{Basic facts about Lagrangian and Hamiltonian systems}
\label{bflt}

In this section we recall some basic facts about Lagrangian and
Hamiltonian dynamics. See for instance \cite{man91}, \cite{bgh98},
\cite{fat06} for detailed proofs. 
Let $M$ be a smooth compact connected
$n$-dimensional manifold without boundary. The elements of the
tangent bundle $TM$ are denoted by $(q,v)$, where $q\in M$ and
$v\in T_q M$. Let $L\in C^2([0,1] \times TM,\R)$ be a Tonelli
Lagrangian, meaning that:

\begin{description}
\item[(L1)] $L$ is $C^2$-strictly convex on the fibers of $TM$, that
  is $\partial_{vv} L >0$;
\item[(L2)] $L$ is superlinear on the fibers of $TM$, that is for any
  $K>0$ there is a finite constant $C(K)$ such that
\[
\forall (t,q,v) \in[0,1] \times TM, \quad L(t,q,v) \geq K |v|_q - C(K)
\]
\end{description}

Here $|v|_q$ is a Riemannian norm on $M$, but by the compactness of 
$M$ condition (L2) does
not depend on the choice of the Riemannian metric on $M$ (up to
changing the constant $C(K)$).

Under assumptions (L1) and (L2), the Legendre transform
\[
\mathscr{L}_L : [0,1] \times TM \rightarrow [0,1]\times T^*M, \quad
(t,q,v) \mapsto (t,q,D_v L(t,q,v)),
\]
is a $C^1$ diffeomorphism. The Fenchel transform of $L$ is the
non-autonomous Hamiltonian on $T^*M$
\[
H (t,q,p) := \max_{v\in T_q M} \left( p[v] - L(t,q,v) \right) =
p[v(t,q,p)] - L(t,q,v(t,q,p)),
\]
where $(t,q,v(t,q,p)) = \mathscr{L}^{-1}_L(t,q,p)$.
Under the above assumptions on $L$, $H$ is a $C^2$ function on
$[0,1]\times T^*M$, it is $C^2$-strictly
convex and superlinear on the fibers of $T^*M$. The associated
non-autonomous Hamiltonian vector field $X_H$ on $T^*M$, defined by
\begin{equation}
\label{hvf}
\omega(X_H(t,q,p),\xi)=-DH(t,q,p)[\xi], \quad \forall (t,q,p)\in
      [0,1]\times T^*M, \; \xi\in T_{(q,p)}T^*M,
\end{equation}
where $\omega=dp\wedge dq$ is the standard symplectic structure on $T^*M$,
is then $C^1$, so it defines a non-autonomous $C^1$ local flow on
$T^*M$. We assume that such a flow is complete:

\begin{description}
\item[(L3)] The solution of
\[
\partial_t \phi^H(t,x)=X_H(t,\phi^H(t,x)), \quad \phi^H(0,x)=x
\]
exists for every $(t,x)\in [0,1]\times T^*M$.
\end{description}

We also use the notation $\phi_t^H(\cdot)=\phi^H(t,\cdot)$.
Assumption (L3) holds, for example, if $H$ satisfies the condition
\begin{equation}
\label{mg}
\partial_t H(t,q,p) \leq c(1+H(t,q,p)) \quad
\forall (t,q,p) \in [0,1] \times T^*M.
\end{equation}\
Indeed, since $DH[X_H]=0$, (\ref{mg}) implies
\[
\frac{d}{dt} H(t,\phi^H(t,x))=\partial_t H(t,\phi^H(t,x)) \leq c(1+H(t,\phi^H(t,x))),
\]
so, by Gronwall Lemma, $H$ is bounded along the flow, which then
exists for every $t\in [0,1]$, by the coercivity of $H$. Condition
(\ref{mg}), written in terms of $L$, becomes
\[
-\partial_t L(t,q,v) \leq c\bigl(1+D_v L(t,q,v)[v] - L(t,q,v)\bigr)
\quad \forall (t,q,v) \in [0,1] \times TM.
\]
Notice that if $L$ - and thus $H$ - is $1$-periodic in time,
condition (L3) is equivalent to the existence of the flow for
every $(t,x)\in \R \times T^*M$. The corresponding flow on $TM$,
obtained by conjugating $\phi^H$ by the Legendre transform
$\mathscr{L}_L$, is denoted by
\[
\phi^L : [0,1] \times TM \rightarrow TM.
\]
Its orbits have the form $t\mapsto (\gamma(t),\gamma'(t))$, where
$\gamma\in C^2([0,1],M)$ solves the Euler-Lagrange equation, which in
local coordinates is written as
\begin{equation}
\label{lageq}
\frac{d}{dt} \partial_v L(t,\gamma(t),\gamma'(t)) = \partial_q
L(t,\gamma(t),\gamma'(t)).
\end{equation}
These orbits are precisely the extremal curves of the Lagrangian
action functional
\[
\mathbb{A}_L(\gamma) := \int_0^1 L(t,\gamma(t),\gamma'(t))\, dt.
\]
Assumptions (L1), (L2), and (L3) imply that solutions with bounded
action are $C^2$-bounded:

\begin{lem}
\label{boundspeedL}
Assume that the $C^2$ Lagrangian $L:[0,1]\times TM \rightarrow \R$
satisfies (L1), (L2), (L3). For every $A\in\R$, the set
\[
\set{\gamma\in C^2([0,1],M)}{ \gamma \mbox{ is a solution of
    (\ref{lageq}) with } \mathbb{A}_L(\gamma) \leq A}
\]
is bounded in $C^2([0,1],M)$.
\end{lem}

\begin{proof}
Since
\[
\int_0^1 L(t,\gamma(t),\gamma'(t))\,dt\leq A,
\]
there exists $t_0 \in [0,1]$ such that
\[
L(t_0,\gamma(t_0),\gamma'(t_0))\leq A.
\]
By assumption (L2), we have in particular $L(t,q,v) \geq | v |_q - C(1)$.
So we get
\[
|\gamma'(t_0))|_{\gamma(t_0)} \leq A + C(1).
\]
Since the Lagrangian flow is a globally defined continuous family of
homeomorphisms, the set
\[
K :=\set{\phi_t^L \circ (\phi_s^L)^{-1} (q,v)}{(s,t)\in [0,1]^2, \;
(q,v) \in TM, \; |v|_q \leq A + C(1)}
\]
is compact. The point $(\gamma(t),\gamma'(t))=\phi^L_{t} \circ
(\phi_{t_0}^L)^{-1}(\gamma(t_0),\gamma'(t_0))$ belongs to $K$ for
every $t\in [0,1]$, so we have an uniform bound in $C^1$. By
computing the time derivative in left-hand side of (\ref{lageq})
and by using assumption (L1), we get in local coordinates
\[
\gamma''(t)=(\partial_{vv}L)^{-1}
\Bigl[ \partial_q L- \partial_t \partial_v L- \partial_{qv} L \,
[\gamma'(t)]\Bigr],
\]
everything being evaluated at $(t,\gamma(t),\gamma'(t))$.
This gives a uniform bound in $C^2$.
\end{proof} \qed

Let $Q$ be a smooth closed submanifold of $M\times M$. We are
interested in the solutions $\gamma:[0,1]\rightarrow M$ of
(\ref{lageq}) which satisfy the nonlocal boundary condition
\begin{equation}
\label{bdry}
\begin{split}
(\gamma(0),\gamma(1)) & \in Q, \\
D_v L(0,\gamma(0),\gamma'(0))[\xi_0] & = D_v
L(1,\gamma(1),\gamma'(1))[\xi_1], \quad \forall (\xi_0,\xi_1) \in
T_{(\gamma(0),\gamma(1))} Q. \end{split}
\end{equation}
Equivalently, the corresponding orbit $x:[0,1]\rightarrow T^*M$ of
$X_H$, that is $(t,x(t))=\mathscr{L}_L(t,\gamma(t),\gamma'(t))$, is
required to satisfy the nonlocal Lagrangian condition
\[
(x(0),-x(1)) \in N^* Q,
\]
where $N^* Q\subset T^*(M \times M)$ denotes the conormal bundle of $Q$ (we
recall that the conormal bundle of a submanifold $Q$ of a manifold $N$
- here $N=M \times M$ - is the Lagrangian submanifold of $T^*N$
consisting of all $(q,p)$ with $q\in Q$ and $p$ vanishing on $T_q Q$).

For instance, if $Q=Q_0 \times Q_1$ is the product of two
submanifolds of $M$, we are looking at solutions $\gamma$ which
join $Q_0$ to $Q_1$, and are normal to $Q_0$ and $Q_1$ in a sense
specified by the Lagrangian $L$. In particular, if $Q_0$ and $Q_1$
are two points, we are looking at solutions joining them without
any further condition. If $L$ is assumed to be $1$-periodic in the
time variable $t$, and $Q=\Delta$ is the diagonal submanifold of
$M\times M$, the solutions of (\ref{lageq}), (\ref{bdry}) are
exactly the $1$-periodic solutions of the Lagrange equation.
Finally, if $Q=M\times M$, condition (\ref{bdry}) reduces to the
Neumann condition
\begin{equation}
\label{neumann}
D_v L(0,\gamma(0),\gamma'(0))= 0, \quad D_v
L(1,\gamma(1),\gamma'(1)) = 0,
\end{equation}
or equivalently
\[
L(0,\gamma(0),\gamma'(0)) = \min_{v\in T_{\gamma(0)} M}
L(0,\gamma(0),v), \quad
L(1,\gamma(1),\gamma'(1)) = \min_{v\in T_{\gamma(1)} M}
L(1,\gamma(1),v).
\]

In order to write the first variation of the action at a curve
$\gamma\in C^1([0,1],M)$, it is convenient to consider an open subset
$U$ of $\R^n$ and a smooth local coordinate system
\begin{equation}
\label{chart}
[0,1] \times U \rightarrow [0,1] \times M, \quad (t,q) \mapsto
(t,\varphi(t,q)),
\end{equation}
such that $\gamma(t) = \varphi(t,\tilde{\gamma}(t))$ for every $t\in
[0,1]$, for some $\tilde{\gamma} \in C^1([0,1],U)$ (see the appendix
for a possible construction of such a map).
Its differential, that is the tangent bundle coordinate system
\[
\Phi: [0,1] \times U \times \R^n \rightarrow [0,1] \times TM, \quad
\Phi(t,q,v) = (t,\varphi(t,q),D\varphi(t,q)[1,v]),
\]
allows to pull back the Lagrangian $L$ onto $[0,1]\times U \times
\R^n$, by setting $\tilde{L} = L \circ \Phi$. If $\xi$ is a variation
of $\gamma$, that is a $C^1$ section of the vector bundle
$\gamma^*(TM)$ over $[0,1]$, and $\tilde{\xi}(t) = D_q \varphi(t,
\tilde{\gamma}(t))^{-1} [\xi(t)]$ is the corresponding variation of
$\tilde{\gamma}$, the first variation of $\mathbb{A}_L$ at $\gamma$
along $\xi$ is
\begin{equation}
\label{firstvar}
\begin{split}
d\mathbb{A}_L(\gamma)[\xi] = d\mathbb{A}_{\tilde L}(\tilde \gamma)
[\tilde \xi] =&\int_0^1  \left( D_q \tilde
L(t,\tilde\gamma,\tilde\gamma')[\tilde\xi] +
D_v \tilde L(t,\tilde\gamma,\tilde\gamma')[\tilde\xi'] \right)  \, dt \\
=& \int_0^1 \left( D_q \tilde L(t,\tilde\gamma,\tilde\gamma') -
  \frac{d}{dt} D_v \tilde L(t,\tilde\gamma,\tilde\gamma')
  \right)[\tilde\xi] \, dt \\
&+ D_v L(1,\gamma(1),\gamma'(1))[\xi(1)] -
D_v L(0,\gamma(0),\gamma'(0))[\xi(0)],
\end{split}
\end{equation}
the last identity holding if $\gamma$, and thus also
$\tilde\gamma$, is $C^2$. Together with a standard regularity
argument, this formula shows that the solutions of (\ref{lageq}),
(\ref{bdry}) are precisely the extremal curves of $\mathbb{A}_L$
on the space of curves $\gamma\in C^1([0,1],M)$ such that
$(\gamma(0),\gamma(1)) \in Q$.

Let $\gamma$ be a solution of (\ref{lageq}), (\ref{bdry}). The
second variation of $\mathbb{A}_L$ at $\gamma$ is given by the
formula:
\begin{equation}
\label{secondvar}
\begin{split}
d^2 \mathbb{A}_L (\gamma)[\xi,\eta] = d^2 \mathbb{A}_{\tilde{L}}
  (\tilde\gamma)[\tilde\xi,\tilde\eta] =
  \int_0^1 \Bigl( D_{vv} \tilde L(t,\tilde\gamma,\tilde\gamma')
[\tilde\xi',\tilde \eta'] + D_{qv} \tilde
  L(t,\tilde\gamma,\tilde\gamma')
[\tilde\xi,\tilde \eta'] \\
+ D_{vq} \tilde
  L(t,\tilde\gamma,\tilde\gamma')
[\tilde\xi',\tilde \eta] +
D_{qq} \tilde L(t,\tilde\gamma,\tilde\gamma') [\tilde\xi,\tilde \eta]
  \Bigr) \, dt. \end{split}
\end{equation}
It is a continuous symmetric bilinear form on the Hilbert space
$W$ consisting of the $W^{1,2}$ sections $\xi$ of $\gamma^*(TM)$
such that $(\xi(0),\xi(1)) \in T_{(\gamma(0),\gamma(1))} Q$. This
bilinear form is a compact perturbation of the form
\[
(\xi,\eta) \mapsto
\int_0^1 \left( D_{vv} \tilde L(t,\tilde\gamma,\tilde\gamma')
[\tilde\xi',\tilde \eta'] + \tilde\xi \cdot \tilde \eta \right) \, dt,
\]
which is coercive on $W$ by (L1). Therefore, the Morse index
$m_Q(\gamma,L)$ (resp.\ the large Morse index $m_Q^*(\gamma,L)$),
i.e. the sum of the multiplicities of the negative (resp.\
non-positive) eigenvalues of the selfadjoint operator on $W$
representing $d^2 \mathbb{A}_L(\gamma)$, is finite. Equivalently,
$m_Q(\gamma,L)$ (resp.\ $m_Q^*(\gamma,L)$) is the dimension of a
maximal linear subspace of $W$ on which $d^2\mathbb{A}_L(\gamma)$
is negative (resp.\ non-positive). Moreover, since the elements in
the kernel of $d^2 \mathbb{A}_L(\gamma)$ solve a system of $n$
second order ODEs, we have the bound
\[
0 \leq m_Q^*(\gamma,L) -m_Q(\gamma,L) \leq 2n.
\]

\begin{lem}
\label{boundmorseindex} Assume that the $C^2$ Lagrangian
$L:[0,1]\times TM \rightarrow \R$ satisfies (L1), (L2), (L3). For
every $A\in\R$, there exists $N \in \N$ such that, for every
$\gamma$ solution of (\ref{lageq}), (\ref{bdry}) with $\mathbb
A_L(\gamma) \leq A$, there holds $m_Q^*(\gamma,L) \leq N$.
\end{lem}

\begin{proof}
By Lemma \ref{boundspeedL} and Ascoli-Arzel\`a theorem, the set of
solutions of (\ref{lageq}), (\ref{bdry}) with $\mathbb A_L(\gamma)
\leq A$ is compact in the $C^1$ topology. Formula
(\ref{secondvar}) shows that the map $\gamma \mapsto d^2
\mathbb{A}_L(\gamma)$ is continuous from the $C^1$ topology of
curves to the operator norm topology of $W^{1,2}$. Since the large
Morse index is upper semi-continuous in the operator norm
topology, the thesis follows.
\end{proof} \qed

\section{Statement of the Lagrangian results}
\label{lagstat} Let $C_Q([0,1],M)$ be the space of continuous
paths $\gamma:[0,1]\rightarrow M$ such that
$(\gamma(0),\gamma(1))\in Q$, endowed with the $C^0$ topology. We
recall that the cuplength of a topological space $X$ is the number
\begin{eqnarray*}
\mathrm{cuplength}\, (X) = \sup \{m\in \N\, | \, \exists \,
  \omega_1,\dots,\omega_m \in H^*(X), \; \deg \omega_j \geq 1 \;
  \forall j=1,\dots,m, \\
\mbox{such that } \omega_1 \cup \dots \cup \omega_m \neq 0\}.
\end{eqnarray*}
Here $H^*(X)$ is the singular cohomology ring of $X$ with integer
coefficients, and $\cup$ denotes the cup product. The main result
of this note concerning Lagrangian systems is the following:

\renewcommand{\theenumi}{\alph{enumi}}
\renewcommand{\labelenumi}{(\theenumi)}

\begin{thm}
\label{main}
Let $M$ be a compact manifold, $Q$ a smooth closed submanifold of
$M\times M$, and $L\in C^2([0,1]\times TM,\R)$ a Lagrangian satisfying (L1),
(L2), and (L3).
\begin{enumerate}
\item The Lagrangian non-local boundary value problem
(\ref{lageq}), (\ref{bdry}) has at least
\[
\mathrm{cuplength}\, (C_Q([0,1],M)) + 1
\]
many solutions. \item If the $k$-th singular homology group
$H_k(C_Q([0,1],M))$ is non-trivial, then problem (\ref{lageq}),
(\ref{bdry}) has a solution $\gamma$ with Morse index estimates
\[
m_Q(\gamma,L) \leq k \leq m_Q^*(\gamma,L).
\]
\item If the $k$-th singular group $H_k(C_Q([0,1],M))$ is non-trivial
  for infinitely many natural numbers $k$,  then problem
(\ref{lageq}), (\ref{bdry}) has an infinite sequence
  of solutions with diverging action and diverging Morse index.
\end{enumerate}
\end{thm}

In the case $Q=\Delta$, the diagonal in $M \times M$, $C_{\Delta}(
[0,1],M)$ coincides with $\Lambda(M)$, the space of free loops in
$M$. If $M$ is compact and simply connected, $H_k(\Lambda(M))\neq 0$
for infinitely many natural numbers $k$ (as proved by Sullivan in
\cite{sul75}), so conclusion (c) of the
above theorem holds. If $M$ is a compact manifold with finite
fundamental group, the above considerations can be applied to its
universal covering (which is still compact), producing contractible periodic
solutions on $M$. These considerations yield to the following:

\begin{cor}
\label{cor1}
Assume that the compact manifold $M$ has finite fundamental group, and let
$L:\R \times TM \rightarrow \R$ be a $C^2$ Lagrangian $1$-periodic in
time and satisfying (L1), (L2), and (L3). Then the Lagrange equation
(\ref{lageq}) has an infinite sequence of $1$-periodic contractible
solutions, with diverging action and diverging Morse index.
\end{cor}

We remark that, in the autonomous case, the corollary above still gives
infinitely many solutions, which are distinct as curves in the phase space
$TM$ since their action diverges. However, in the case of the geodesic flow,
all these solutions could be the same geodesic parametrized with
increasing speed.

When $M$ is simply connected, also $\Omega(M)$ - the space of based
loops in $M$ - has infinitely many non-vanishing homology groups (a
classical result by Serre, \cite{ser51}).
Moreover, if $Q$ is the singleton $\{(q_0,q_1)\}$, with
$q_0,q_1\in M$, the space $C_Q([0,1],M)$ is homotopically
equivalent to $\Omega(M)$. Therefore,
conclusion (c) of Theorem \ref{main} implies the following:

\begin{cor}
\label{cor2}
Assume that the compact manifold $M$ has finite fundamental group, and let
$L:[0,1] \times TM \rightarrow \R$ be a $C^2$ Lagrangian satisfying
(L1), (L2), and (L3). Then for every pair of points $q_0,q_1$ in $M$
and for every continuous path $\overline{\gamma}$ joining them, there is
an infinite sequence of solutions $\gamma$ of (\ref{lageq}) joining
$q_0$ and $q_1$, homotopic to $\overline{\gamma}$, with diverging
action and diverging Morse index.
\end{cor}

In the case $Q=Q_0 \times M$, with $Q_0$ a closed submanifold of $M$,
the map
\[
C_{Q_0 \times M} ([0,1],M) \rightarrow Q_0, \quad \gamma \mapsto
\gamma(0),
\]
is a homotopy equivalence, a homotopy inverse of it being the
function mapping each $q\in Q_0$ into the constant path $\gamma(t)\equiv q$.
Therefore statement (a) in Theorem \ref{main} has the following
consequence:

\begin{cor}
\label{cor3}
Assume that $M$ is a compact manifold, and let
$L:[0,1] \times TM \rightarrow \R$ be a $C^2$ Lagrangian satisfying
(L1), (L2), and (L3). Then the Lagrange equation (\ref{lageq}) has at
least $\mathrm{cuplength}\, (Q_0)+1$ many solutions $\gamma$ satisfying the
boundary conditions
\[
\gamma(0)\in Q_0, \quad D_v L(0,\gamma(0),\gamma'(0))|_{T_{\gamma(0)}
  Q_0} = 0, \quad
D_v L(1,\gamma(1),\gamma'(1)) = 0.
\]
In particular, the Neumann problem (\ref{neumann}) for the
Lagrange equation (\ref{lageq}) has at least
$\mathrm{cuplength}\,(M)+1$ many solutions.
\end{cor}

\section{Quadratic Lagrangians}
\label{qls}

Let us assume that the $C^2$ Lagrangian $L:[0,1] \times M \rightarrow
\R$ satisfies the conditions:
\begin{description}
\item[(L1')] There is a constant $\ell_0>0$ such that
$\partial_{vv} L(t,q,v)
  \geq \ell_0 I$.
\item[(L2')] There is a constant $\ell_1>0$ such that
\[
\|\partial_{vv} L(t,q,v)\| \leq \ell_1, \quad \|\partial_{qv} L(t,q,v)\| \leq
\ell_1 (1+|v|_q),
\quad \|\partial_{qq} L(t,q,v)\| \leq \ell_1(1+|v|_q^2).
\]
\end{description}
Condition (L1') is expressed in term of a Riemannian metric on
$M$, but the condition does not depend on such a metric (up to
changing the constant $\ell_0$). Condition (L2') is expressed in
terms of a system of local coordinates on $M$ and of a Riemannian
metric on $M$, but the condition does not depend on these choices
(up to changing the constant $\ell_1$). Assumption (L1') implies
that $L$ grows at least quadratically in $v$, while (L2') implies
that $L$ grows at most quadratically.

We observe that assumption (L1') implies both (L1) and (L2). Under
these assumptions, Benci \cite{ben86} has shown that the $W^{1,2}$
functional setting used in the study of geodesics allows to prove
existence result for the critical points of the action functional
$\mathbb{A}_L$. The aim of this section is to recall Benci's
results, extending them from the periodic case to the case of a
general boundary condition of the form (\ref{bdry}).

Let $W^{1,2}([0,1],M)$ be the space of all curves $\gamma:[0,1]
\rightarrow M$ of Sobolev class $W^{1,2}$. This space has a
natural structure of a Hilbert manifold modeled on
$W^{1,2}([0,1],\R^n)$. Indeed, a smooth atlas on
$W^{1,2}([0,1],M)$ is defined by composition by diffeomorphisms of
the form (\ref{chart}). Its tangent space at $\gamma\in
W^{1,2}([0,1],M)$ is identified with the space of $W^{1,2}$
sections of the vector bundle $\gamma^*(TM)$. For sake of
completeness, the following proposition, as well as Proposition
\ref{ps} below, are proved in the appendix.

\begin{prop}
\label{reg}
If $L$ satisfies (L2')
then $\mathbb{A}_L$ is of class $C^2$
on $W^{1,2}([0,1],M)$, and $D \mathbb{A}_L(\gamma)$ coincides with
$d\mathbb{A}_L(\gamma)$, given by (\ref{firstvar}), while
$D^2 \mathbb{A}_L(\gamma)$ at a critical point $\gamma$
coincides with $d^2 \mathbb{A}_L(\gamma)$, given by (\ref{secondvar}).
\end{prop}

Let $W^{1,2}_Q([0,1],M)$ be the set of $\gamma\in W^{1,2}([0,1],M)$
such that $(\gamma(0),\gamma(1))\in Q$. Being the inverse image of $Q$
by the smooth submersion
\[
W^{1,2}([0,1],M) \rightarrow M \times M, \quad \gamma \mapsto
(\gamma(0),\gamma(1)),
\]
$W^{1,2}_Q([0,1],M)$ is a closed smooth submanifold of
$W^{1,2}([0,1],M)$. By Proposition \ref{reg} and identity
(\ref{firstvar}), the restriction of $\mathbb{A}_L$ is $C^2$, and
the critical points of such a restriction are exactly the
solutions of (\ref{lageq}), (\ref{bdry}). The second differential
of $\mathbb{A}_L|_{W^{1,2}_Q([0,1],M)}$ at a critical point
$\gamma$ is the restriction of $D^2 \mathbb{A}_L(\gamma)$ to
$T_{\gamma} W^{1,2}_Q([0,1],M)$, and it coincides with the
corresponding restriction of $d^2 \mathbb{A}_L(\gamma)$, given by
(\ref{secondvar}). In particular the Hessian of
$\mathbb{A}_L|_{W^{1,2}_Q([0,1],M)}$ at a critical point is a
Fredholm operator. Therefore, the Morse index and the large Morse
index of such a critical point are the numbers $m_Q(\gamma,L)$ and
$m_Q^*(\gamma,L)$ defined in section \ref{bflt}.

A Riemannian metric $\langle \cdot,\cdot \rangle$ on $M$ induces a
complete Riemannian metric on the Hilbert manifold $W^{1,2}([0,1],M)$,
namely
\begin{equation}
\label{metric}
\langle \xi,\eta \rangle_{W^{1,2}} := \int_0^1 \left(\langle \nabla_t \xi,
\nabla_t \eta \rangle_{\gamma(t)} + \langle \xi,\eta
\rangle_{\gamma(t)} \right)\, dt , \quad
\forall \gamma\in W^{1,2}([0,1],M), \; \xi,\eta \in T_{\gamma}
W^{1,2}([0,1],M),
\end{equation}
where $\nabla_t$ denotes the Levi-Civita covariant derivative along
$\gamma$. We recall that a $C^1$ functional $f$ on a Riemannian Hilbert
manifold $(\mathscr{M},\|\cdot\|)$ satisfies the Palais-Smale
condition if every sequence $(u_h)\in \mathscr{M}$ such that $f(u_h)$
is bounded and $\|Df(u_h)\|_*$ is infinitesimal is compact (here
$\|\cdot\|_*$ denotes the dual norm on $T_{u_h}^*\mathscr{M}$).

\begin{prop}
\label{ps}
If $L$ satisfies (L1') and (L2'), then the restriction of
$\mathbb{A}_L$ to $W^{1,2}_Q([0,1],M)$ satisfies the Palais-Smale
condition.
\end{prop}

Let $\alpha,\beta\in H_*(X)$ be non-zero singular homology classes in the
topological space $X$, with $\deg \alpha < \deg \beta$.
We recall that $\alpha$ is said to be subordinate to $\beta$, $\alpha
<\beta$, if there exists a singular cohomology class $\omega\in
H^*(X)$ such that $\alpha = \beta \cap \omega$, where $\cap :
H_{p+q}(X) \otimes H^q(X) \rightarrow H_p(X)$ denotes the cap
product. If $\alpha<\beta$ and $b$ is a singular cycle representing
$\beta$, there exists a singular cycle $a$ representing $\alpha$
with support contained in the support of $b$.

The fact that the Lagrangian action functional $\mathbb{A}_L$ is
$C^2$, bounded from below (because so is $L$), satisfies the
Palais-Smale condition on a complete Riemannian manifold,
and has a Fredholm Hessian, implies the
following result:

\begin{thm}
\label{mainq}
Let $M$ be a compact manifold and let $Q$ be a smooth closed submanifold of
$M \times M$. Assume that the $C^2$ Lagrangian
$L:[0,1]\times TM \rightarrow \R$ satisfies (L1') and (L2').
\begin{enumerate}
\item Let $\alpha\in H_k(W^{1,2}_Q([0,1],M))$, $k\in \N$, be a non-zero
homology class, and let $\mathscr{K}_{\alpha}$ be the
family consisting of the supports of the singular cycles in
$W^{1,2}_Q([0,1],M)$ representing $\alpha$. Then the number
\[
c_{\alpha}(L) := \inf_{K\in \mathscr{K}_{\alpha}} \max_{\gamma\in K}
\mathbb{A}_L (\gamma)
\]
is a critical level of the action function $\mathbb{A}_L$, and there
is a critical point $\gamma$ with
\[
\mathbb{A}_L(\gamma) = c_{\alpha}(L), \quad
m_Q(\gamma,L) \leq k \leq m_Q^*(\gamma,L).
\]
\item If $\alpha_1 < \dots < \alpha_m$ are non-zero homology
  classes in $H_*(W^{1,2}_Q([0,1],M))$, then
\[
c_{\alpha_1}(L) \leq \dots \leq c_{\alpha_m} (L),
\]
and either all the inequalities are strict, or one
of the critical levels $c_{\alpha_j}(L)$, $j=2,\dots,m$, contains a
continuum of critical points.
\end{enumerate}
\end{thm}

Indeed, everything follows from the abstract results of section I.3.2
in \cite{cha93}, apart from the Morse index estimate which is proved
in \cite{vit88}.

\section{Convex quadratic modifications}

Let $L\in C^2([0,1]\times TM,\R)$ be a Lagrangian satisfying (L1),
(L2), and (L3).

\begin{defn}
\label{qrm}
We say that a Lagrangian $L_0\in C^2([0,1]\times TM, \R)$ is a convex
quadratic $R$-modification of $L$ if:
\begin{enumerate}
\item $L_0(t,q,v) =  L(t,q,v)$ for $|v|_q \leq R$;
\item $L_0$ satisfies (L1') and (L2');
\item $L_0(t,q,v) \geq |v|_q - C(1)$, where $C(1)$ is defined in
condition (L2).
\end{enumerate}
\end{defn}

Constructing a convex quadratic $R$-modifications of a Tonelli
Lagrangian $L$ is not difficult, although some care is needed in
order to preserve the convexity and the linear lower bound. Here is a
possible construction.

Let $\varphi:\R \rightarrow \R$ be a smooth function such that
$\varphi(s)=s$ for $s\leq 1$ and $\varphi(s)$ is constant for $s\geq
2$. We fix a positive number $\lambda$ such that
\begin{equation}
\label{eq1}
\lambda \geq \max \set{L(t,q,v)}{(t,q,v)\in [0,1]\times TM, \; |v|_q
\leq 2R},
\end{equation}
and we define a new Lagrangian $L_1$ on $[0,1]\times TM$ by
\[
L_1 := \lambda \varphi \left( \frac{L}{\lambda} \right).
\]
Since $L$ is coercive, $L_1$ is constant outside a compact set, so we
can find a positive number $\mu$ such that
\begin{equation}
\label{eq2}
\mu I \geq - \partial_{vv} L_1, \quad \mbox{on } [0,1] \times TM.
\end{equation}
Up to replacing $\mu$ with a larger number, we may also assume that
\begin{equation}
\label{eq3}
4R\mu \geq 1, \quad 2 R^2 \mu \geq 2R - C(1) - \min L_1.
\end{equation}
The affine function $s\mapsto \mu s - 2 \mu R^2$ is negative on
$(-\infty,R^2]$ and positive on $[4R^2,+\infty)$, so we can find a 
smooth convex function $\psi:\R
\rightarrow \R$ such that $\psi(s)=0$ for $s\leq R^2$ and $\psi(s) =
\mu s - 2 \mu R^2$ for $s\geq 4R^2$. We define the Lagrangian $L_0$
by
\[
L_0 (t,q,v) := L_1(t,q,v) + \psi\left(|v|_q^2 \right).
\]
By (\ref{eq1}), 
\begin{equation}
\label{eq4}
L_0(t,q,v) = L(t,q,v) + \psi\left(|v|_q^2 \right), \quad \mbox{for }
|v|_q \leq 2R,
\end{equation}
and in particular $L_0$ coincides with $L$ for $|v|_q\leq R$, proving
property (a) of Definition \ref{qrm}. By (\ref{eq4}) and the convexity
and the monotonicity of $\psi$, 
\begin{equation}
\label{eq5}
\partial_{vv} L_0 (t,q,v) \geq \partial_{vv} L(t,q,v), \quad \mbox{for
} |v|_q \leq 2R.
\end{equation}
On the other hand, by (\ref{eq2}),
\begin{equation}
\label{eq6}
\partial_{vv} L_0 (t,q,v) = \partial_{vv} L_1(t,q,v) + 2 \mu I \geq \mu
I, \quad \mbox{for } |v|_q \geq 2R.
\end{equation}
Since $L$ is $C^2$-strictly convex, (\ref{eq5}) and (\ref{eq6}) imply
that $L_0$ satisfies (L1'). Since $L_1$ is constant outside a compact set,
\[
L_0(t,q,v) = \mu |v|_q^2 + \mbox{constant},
\]
for $|v|_q$ large, so $L_0$ satisfies also (L2'), concluding the
proof of property (b) in Definition \ref{qrm}. By (\ref{eq4}) and
(L2),
\begin{equation}
\label{eq7}
L_0(t,q,v) \geq L(t,q,v) \geq |v|_q - C(1), \quad \mbox{for } |v|_q
\leq 2R.
\end{equation}
On the other hand, by (\ref{eq3}),
\begin{equation}
\label{eq8}
L_0(t,q,v) \geq \min L_1 + \mu |v|_q^2 - 2 \mu R^2 \geq 
|v|_q - C(1), \quad \mbox{for } |v|_q \geq 2R.
\end{equation}
Indeed, (\ref{eq3}) implies that
\[
\mu s^2 - s + \min L_1 - 2 \mu R^2 + C(1) \geq 0, \quad
\mbox{for } s\geq 2R,
\]
as it easily seen by evaluating the above polynomial and its
derivative at $s=2R$. Inequalities (\ref{eq7}) and (\ref{eq8}) imply
that $L_0$ satisfies property (c) of Definition \ref{qrm}. Therefore,
$L_0$ is a convex quadratic $R$-modification of $L$. 

If $L$ is $1$-periodic
in time, the above construction produces a Lagrangian which is also
$1$-periodic in time.

\begin{lem}
\label{boundspeedqRm}
Let $M$ be a compact manifold and $L\in C^2([0,1]\times TM,\R)$ be a
Lagrangian satisfying (L1), (L2), (L3).
For every $A>0$ there exists a number $R(A)$ such that for any
$R > R(A)$ and  for any Lagrangian $L_0$ which is a convex quadratic
$R$-modification of $L$, the following holds:
if $\gamma$ is a critical point of $\mathbb{A}_{L_0}$ such that
$\mathbb{A}_{L_0}(\gamma) \leq A$, then $\|\gamma'\|_{\infty}
\leq R(A)$. In particular, such a $\gamma$ is an
extremal curve of
$\mathbb{A}_L$, and $\mathbb{A}_{L_0}(\gamma)=\mathbb{A}_{L}(\gamma)$.
\end{lem}

\begin{proof}
Let $C(1)$ be the constant appearing in assumption (L2).
Since the Lagrangian flow of $L$ is globally defined, as in the proof of
Lemma \ref{boundspeedL} we find that the set
\[
K :=\set{\phi_t^L \circ (\phi_s^L)^{-1} (q,v)}{(s,t)\in [0,1]^2, \;
(q,v) \in TM, \; |v|_q \leq A + C(1)}
\]
is compact. Notice that $K$ contains every $(q,v)\in TM$ with
$|v|_q\leq A + C(1)$. Let $R(A)$ be the maximum of $|v|_q$ for
$(q,v)\in K$. Let $R
> R(A)$, $L_0$ a convex quadratic $R$-modification of $L$, and
$\gamma$ a critical point of $\mathbb{A}_{L_0}$ such that
$\mathbb{A}_{L_0}(\gamma) \leq A$. Then there exists $t_0 \in
[0,1]$ such that
\[
|\gamma'(t_0)|_{\gamma(t_0)} - C(1)\leq
L_0(t_0,\gamma(t_0),\gamma'(t_0))\leq \mathbb{A}_{L_0}(\gamma) \leq A,
\]
where the first inequality follows by $(c)$ in Definition \ref{qrm}.
So we get
\[
|\gamma'(t_0))|_{\gamma(t_0)} \leq A + C(1),
\]
which implies that $(\gamma(t_0),\gamma'(t_0)) \in K$.
Let $I \subset [0,1]$ be the maximal interval containing $t_0$ such that
$(\gamma(t),\gamma'(t)) \in K$ for every $t \in I$.
If $t \in I$, then $|\gamma'(t))|_{\gamma(t)} \leq R(A)<R$.
Therefore $L$ and $L_0$ coincide in a neighborhood of
$(t,\gamma(t),\gamma'(t))$, and so do the Lagrangian vector fields
defined by $L$ and $L_0$. Hence, by the definition of $K$,
$(\gamma(s),\gamma'(s))$ belongs to $K$ for $s$ in a neighborhood of
$t$ in $[0,1]$, proving that $I$ is open in $[0,1]$. Being obviously
closed, $I$ coincides with $[0,1]$. Therefore $\|\gamma'\|_{\infty}
\leq R(A)$.
\end{proof} \qed

\section{Proof of Theorem \ref{main}}

We denote by $C^k_Q([0,1],M)$ the space of $C^k$ curves
$\gamma:[0,1]\rightarrow M$ such that $(\gamma(0),\gamma(1))\in Q$,
endowed with the $C^k$ topology, for $1\leq k \leq \infty$.
Let $\alpha\in H_*(C_Q([0,1],M))$ be a non-zero homology class.
Since the inclusions
\begin{equation}
\label{inclsn}
C^k_Q([0,1],M) \hookrightarrow W^{1,2}_Q([0,1],M) \hookrightarrow
C_Q([0,1],M) , \quad 1\leq k \leq \infty,
\end{equation}
are homotopy equivalences, we may regard $\alpha$ as a non-zero
homology class on each of these different spaces of curves. Let
$\mathscr{K}_{\alpha}'$ be the (non-empty) family of supports of
singular cycles in $C^1_Q([0,1],M)$ representing the homology class
$\alpha$. Since singular chains have compact support, and since
$\mathbb{A}_L$ is bounded from below, and continuous on $C^1_Q([0,1],M)$,
\[
c_{\alpha}'(L) := \inf_{K \in \mathscr{K}_{\alpha}'} \max_{\gamma\in K}
\mathbb{A}_L(\gamma)
\]
is a finite number.

Let $m\geq 1$ be an integer, and let $\alpha_1 < \dots < \alpha_m$ be
non-zero elements of $H_*(C_Q([0,1],M))$. Let $K_j\in
\mathscr{K}_{\alpha_j}'$ be such that
\[
\max_{\gamma \in K_j} \mathbb{A}_L(\gamma) \leq c_{\alpha_j}'(L)+1, \quad
\forall j=1,\dots,m.
\]
Set
\[
A:=\max \set{c_{\alpha_j}'(L)+1}{j=1,\dots,m} = c_{\alpha_m}'(L)+1.
\]
We choose $R>0$ so large that $R> R(A)$, the number provided by Lemma
\ref{boundspeedqRm}, and
\begin{equation}
\label{Rg}
R\geq \sup \Bigl\{\|\gamma'\|_{\infty} \, \Big| \,
\gamma\in \bigcup_{j=1}^m K_j\Bigr\}.
\end{equation}
Let $L_0$ be a convex quadratic $R$-modification of $L$.
Since the first inclusion in (\ref{inclsn}) is a homotopy equivalence,
the family $\mathscr{K}_{\alpha_j}'$
is contained in the family $\mathscr{K}_{\alpha_j}$ appearing in Theorem
\ref{mainq}. Therefore, taking also (\ref{Rg}) into account and the
fact that $L_0(t,q,v)$ coincides with $L(t,q,v)$ if $|v|_q\leq R$,
\[
c_{\alpha_j}(L_0)  := \inf_{K \in \mathscr{K}_{\alpha_j}}
\max_{\gamma\in K} \mathbb{A}_{L_0}(\gamma) \leq \max_{\gamma \in
K_j} \mathbb{A}_{L_0} (\gamma) = \max_{\gamma\in K_j} \mathbb{A}_L
(\gamma) \leq c_{\alpha_j}'(L) + 1 \leq A.
\]
By Theorem \ref{mainq} (a) the functional
$\mathbb{A}_{L_0}|_{W^{1,2}_Q([0,1],M)}$ has a critical point
$\gamma_j$ with $\mathbb{A}_{L_0}(\gamma_j) = c_{\alpha_j}(L_0)$ and
\[
m_Q(\gamma_j;L_0) \leq \deg \alpha_j \leq
m^*_Q(\gamma_j;L_0).
\]
Since $\mathbb{A}_{L_0}(\gamma_j) = c_{\alpha_j}(L_0) \leq A$ and
$R> R(A)$, by Lemma \ref{boundspeedqRm} $\|\gamma_j'\|_{\infty}
\leq R(A) < R$. Since $L_0(t,q,v)$ coincides with $L(t,q,v)$ if
$|v|_q\leq R$, the curve $\gamma_j$ is a solution of
(\ref{lageq}), (\ref{bdry}) for the Lagrangian $L$, and
\[
m_Q(\gamma_j;L) = m_Q(\gamma_j;L_0) \leq \deg \alpha_j \leq
m^*_Q(\gamma_j;L_0) = m^*_Q(\gamma_j;L).
\]
By Theorem \ref{mainq} (b), either the critical levels
$c_{\alpha_j}(L_0)$ are all distinct - hence the curves $\gamma_j$
are also distinct - or one of these levels contains a continuum of
critical points of $\mathbb{A}_{L_0}$, and thus of $\mathbb{A}_L$. In any case,
we deduce that the original problem (\ref{lageq}), (\ref{bdry})
has at least $m$ solutions. Statement (a) of Theorem \ref{main}
follows from the identity
\[
\mathrm{cuplength}\, (X) + 1 =
\sup\set{m\in \N}{\exists\, \alpha_1 < \dots < \alpha_m \mbox{ non-zero
    elements of } H_*(X)},
\]
(see \cite{cha93}, Theorem I.1.1).

Taking $m=1$ and $\alpha_1$ a non-zero element of
$H_k(C_Q([0,1],M))$ in the above argument proves statement (b) of
Theorem \ref{main}. Under the assumption of Theorem \ref{main}
(c), problem (\ref{lageq}), (\ref{bdry}) has then a sequence of
solutions with diverging large Morse index. By Lemma
\ref{boundmorseindex}, the fact that the large Morse index
diverges implies that also the action diverges. This concludes the
proof of Theorem \ref{main}.

\begin{rem}
A natural question is whether $c'_{\alpha_j}(L)$ is a critical
value of $\mathbb{A}_L$ or not. We observe that, choosing compact
sets $K_j\in \mathscr{K}_{\alpha_j}'$ such that
\[
\max_{\gamma \in K_j} \mathbb{A}_L(\gamma) \leq c_{\alpha_j}'(L)+\epsilon
\quad\forall j=1,\dots,m
\]
for a certain $\epsilon>0$ small, by the above argument we get
\[
\mathbb{A}_{L_\epsilon}(\gamma_\epsilon)=c_{\alpha_j}(L_\epsilon)
\leq c_{\alpha_j}'(L) + \epsilon,
\]
where $L_\epsilon$ is a $R_\epsilon$-modification of $L$ and
$\gamma_\epsilon$ is a critical point for both $L$ and
$L_\epsilon$. Now, since $\mathbb{A}_L(\gamma_\epsilon)=\mathbb{A}_{L_\epsilon}(\gamma_\epsilon)\leq
c_{\alpha_j}'(L) + \epsilon$, the family of critical points
$(\gamma_\epsilon)$ is bounded in $C^2$ (by Lemma
\ref{boundspeedL}), so, up to subsequences, it converges in the
$C^1$ topology to a critical point $\gamma_\infty$ which satisfies
$\mathbb{A}_L(\gamma_\infty) \leq c_{\alpha_j}'(L)$.

Let us now assume that $L$ has at most quadratic growth in $v$,
that is $L(t,q,v) \leq C(1+|v|_q^2)$ for a certain constant $C$.
In this case, up to enlarging $R$, we can assume that every
$R$-modification $L_0$ of $L$ satisfies $L_0 \geq L$. This implies
that $\mathbb A_{L_0} \geq \mathbb A_L$, so we get
\[
c_{\alpha_j}(L_0)  = \inf_{K \in \mathscr{K}_{\alpha_j}}
\max_{\gamma\in K} \mathbb{A}_{L_0}(\gamma) = \inf_{K \in
\mathscr{K}_{\alpha_j}'} \max_{\gamma\in K}
\mathbb{A}_{L_0}(\gamma) \geq \inf_{K \in \mathscr{K}_{\alpha_j}'}
\max_{\gamma\in K} \mathbb{A}_{L}(\gamma)=c_{\alpha_j}'(L),
\]
where in the second equality we used the density of $C^1_Q([0,1],M)$
in $W^{1,2}_Q([0,1],M)$
in the $W^{1,2}$ topology, and the continuity of $\mathbb{A}_{L_0}$
with respect to this topology.
So, if $\gamma_{\epsilon_k}$ converges to $\gamma_\infty$ as 
$k \to \infty$, we get in this case
\[
c_{\alpha_j}'(L) \geq \mathbb{A}_L(\gamma_\infty)= \lim_k \mathbb A_L(\gamma_{\epsilon_k})
= \lim_k c_{\alpha_j}(L_{\epsilon_k}) \geq c_{\alpha_j}'(L),
\]
which implies that the number $c_{\alpha_j}'(L)$ is indeed a
critical value for the action of $L$. It is not clear to the
authors whether this remains true for any Tonelli Lagrangian.
\end{rem}

\section{Statement of the Hamiltonian Results}
\label{secHamilt}

Let $H\in C^2([0,1]\times T^*M,\R)$ be a time-dependent Hamiltonian on
the cotangent bundle of the smooth compact connected $n$-dimensional
manifold $M$. We denote by $\lambda=p dq$ the Liouville 1-form on
$T^*M$, whose differential $d\lambda = dp \wedge dq$ is the standard
symplectic form $\omega$ on $T^*M$. The Liouville vector field
$Y$ on $T^* M$ is defined by
\begin{equation}
\label{lvf}
\omega (Y, \cdot) = \lambda,
\end{equation}
which in local coordinates becomes $Y=\sum_i p_i \partial_{p_i}$.
We denote by $X_H$ the Hamiltonian vector field defined by
(\ref{hvf}), and by $\phi^H$ its non-autonomous local flow.

We assume that the Hamiltonian $H$ satisfies the following conditions:
\begin{description}

\item[(H1)] The action integrand $DH[Y]-H$ is coercive on
$[0,1]\times T^*M$, that is
\[
DH(t,q,p)[Y(q,p)]-H(t,q,p) \geq a(|p|_q), \quad \text{where }
\lim_{s\rightarrow +\infty} a(s)=+\infty;
\]

\item[(H2)] The function $H$ is superlinear on the fibers of $T^*M$,
 that is
 \[
H(t,q,p)\geq h(|p|_q), \quad \text{where } \lim_{s\rightarrow
+\infty}\frac{h(s)}{s} = +\infty;
\]

\item[(H3)] The non-autonomous flow $\phi^H$ is globally defined
on $[0,1]\times T^*M$.
\end{description}

\begin{rem}
\label{rmkH12}
For fixed $t \in [0,1]$ and $(q,p)  \in T^*M$ with $|p|_q=1$, let
us consider the function $f:[0,+\infty) \to \R$ defined as
\[
f(s):=DH(t,q,sp)[Y(q,sp)]-H(t,q,sp).
\]
Then, if we set $g(s):=H(t,q,sp)$, we get $g'(s)s - g(s)=f(s)$,
from which
\[
\frac{g(s)}{s} = g(1) + \int_1^s \frac{f(\sigma)}{\sigma^2}
\,d\sigma.
\]
Thus we see that $g(s)$ is superlinear if and only if $f(s)/s^2$
is not integrable at infinity. Therefore, assumption (H2) implies
that $f(s)/s^2$ is not integrable at infinity. Hence (H2) is in a
certain sense ``stronger'' than (H1), and the function $a(s)$ appearing
in (H1) is expected to grow at least linearly.
\end{rem}

Let $Q$ be a smooth closed submanifold of $M\times M$. We are
interested in the solutions $x:[0,1]\rightarrow T^*M$ of
\begin{equation}
\label{he}
x'(t) = X_H(t,x(t)),
\end{equation}
satisfying the boundary condition
\begin{equation}
\label{hbc}
(x(0),-x(1)) \in N^* Q,
\end{equation}
where $N^* Q \subset T^*(M\times M)$ denotes the conormal bundle
of $Q$ (see section \ref{bflt}).

The Hamiltonian action functional acts on paths on $T^* M$, and it is
defined as
\[
\mathbb{A}_H(x) = \int_0^1 \bigl( \lambda[x'(t)] -
H(t,x(t)) \bigr) \, dt,
\]
where $x:[0,1]\rightarrow T^*M$. Notice that if
$x:[0,1]\rightarrow T^* M$ is an orbit of $X_H$ then by
(\ref{lvf}), (\ref{he}) and (\ref{hvf}) we have
\begin{equation}
\label{actint} \lambda [x'] - H(t,x) = \omega(Y(x),X_H(t,x)) -
H(t,x) = DH(t,x)[Y(x)] - H(t,x).
\end{equation}
The above identity explains why we refer to the quantity $DH[Y]-H$
as to the {\em action integrand}. Notice also that if $H$ is the
Fenchel transform of a Tonelli Lagrangian $L$, then the value of the
action integrand along an orbit coincides pointwise with the value
of the Lagrangian along the corresponding orbit on the tangent
bundle.

The first variation of $\mathbb{A}_H$ on the space of free paths
on $T^*M$ is
\begin{equation}
\label{fvhaf}
d\mathbb{A}_H(x)[\xi] = \int_0^1 \bigl( \omega(\xi,x') - D_x
H(t,x)[\xi] \bigr)\, dt + \lambda(x(1))[\xi(1)] -
\lambda(x(0))[\xi(0)],
\end{equation}
where $\xi$ is a section of $x^*(TT^*M)$.
The functional $\mathbb{A}_H$ is unbounded from above and from
below. However, we have the following:

\begin{lem}
\label{hae}
Assume that the $C^2$ Hamiltonian $H:[0,1]\times T^*M \rightarrow \R$
satisfies (H1), (H3). For every $A\in \R$ the set
\[
\set{ x\in C^2([0,1],T^*M)}{x \mbox{ is a solution of (\ref{he}) with
  } \mathbb{A}_H(x) \leq A}
\]
is bounded in $C^2([0,1],T^*M)$.
\end{lem}

\begin{proof}
Since
\[
\int_0^1 \bigl( \lambda[x'(t)] - H(t,x(t)) \bigr) \, dt \leq A,
\]
there exists $t_0\in [0,1]$ such that
\[
DH(t_0,x(t_0))[Y(x(t_0))] - H(t_0,x(t_0)) = \lambda[x'(t_0)] -
H(t_0,x(t_0)) \leq A,
\]
where we have used (\ref{actint}) in the first equality. By
assumption (H1), $x(t_0)=(q(t_0),p(t_0))$ belongs to the compact
set $\{(p,q) \in T^*M \mid a(|p|_q) \leq A\}$. Since the
Hamiltonian flow is a globally defined continuous family of
homeomorphisms, also the set
\[
K := \set{\phi_t^H \circ (\phi_{s}^H)^{-1} (p,q)}{(s,t)\in [0,1]^2, \
  (q,p) \in T^*M, \ a(|p|_q) \leq A}
\]
is compact. The point $x(t)=\phi_t^H\circ (\phi_{t_0}^H)^{-1}
(x(t_0))$ belongs to $K$ for every $t\in [0,1]$, so we have a
uniform bound for $x$ in $C^0$. Since $x$ solves an ordinary
differential equation with $C^1$ coefficients, we deduce a uniform
bound in $C^2$.
\end{proof} \qed

The aim of the remaining part of this note is to prove the following:

\begin{thm}
\label{main2} Let $M$ be a compact manifold, $Q$ a smooth closed
submanifold of $M\times M$, and $H \in C^2([0,1]\times T^*M,\R)$ a
Hamiltonian satisfying (H1), (H2), and (H3).
\begin{enumerate}
\item The non-local boundary value problem (\ref{he}), (\ref{hbc})
has at least
\[
\mathrm{cuplength}_{\Z_2} \, (C_Q([0,1],M)) + 1
\]
many solutions.
\item If the $k$-th singular group $H_k(C_Q([0,1],M),\Z_2)$ is non-trivial
for infinitely many natural numbers $k$,  then problem
(\ref{he}), (\ref{hbc}) has an infinite sequence
of solutions with diverging action.
\end{enumerate}
\end{thm}

The Hamiltonian analogues of Corollaries \ref{cor1}, \ref{cor2},
\ref{cor3} are easily deduced from the above theorem.

\begin{rem}
The appearance of cohomology with $\Z_2$ coefficients in the above
theorem is related to the fact that the proof makes use of degree
theory for Fredholm maps. The extension to integer valued
cohomology could be probably obtained by considering the concept
of determinant bundle over the space of Fredholm operators (see
\cite{fh93}). The statements concerning the Morse index of
solutions in Theorem \ref{main} can also be extended to the
Hamiltonian setting, by replacing the Morse index by the Maslov
index of suitable paths of Lagrangian subspaces (see for instance
\cite{sz92}, \cite{web02}, and \cite{as06}). 
Since the Fenchel transform of a Tonelli Lagrangian
satisfies (H1) and (H2), and since the Hamiltonian Maslov index
coincides with the Lagrangian Morse index in this case, such an
extension of Theorem \ref{main2} implies Theorem \ref{main}.
\end{rem}

\section{Quadratic Hamiltonians}

In \cite{cie94} Cieliebak extended Floer's idea of using
$J$-holomorphic curves in the study of periodic orbits of Hamiltonian
systems on compact manifolds to the case of cotangent bundles. In this
section we describe Cieliebak's results, in the setting of our more
general boundary conditions.

Let us assume that the Hamiltonian $H:[0,1]\times T^*M
\rightarrow \R$ is smooth and satisfies the conditions:
\begin{description}
\item [(H1')] There are constants $h_0>0$, $h_1\geq 0$, such that
\[
DH(t,q,p)[Y(q,p)] - H(t,q,p) \geq h_0 |p|_q^2 - h_1,
\]
for every $(t,q,p)\in T^*M$.
\item [(H2')] There is a constant $h_2\geq 0$ such that
\[
\|\partial_q H(t,q,p)\| \leq h_2 (1+|p|_q^2), \quad
 \|\partial_p H(t,q,p)\| \leq h_2 (1+|p|_q),
\]
for every $(t,q,p)\in T^*M$.
\end{description}

Condition (H1') is expressed in term of a norm on the vector
bundle $T^*M$, but the condition does not depend on such a norm
(up to changing the constants). Condition (H2') uses also a system
of local coordinates on $M$ and the induced system on $T^*M$, but
the condition does not depend on these choices. Condition (H1')
implies that $H$ grows at least quadratically on the fibers 
(see Remark \ref{rmkH12}), while (H2') implies that $H$ grows at most 
quadratically. In
particular, (H1') implies both (H1) and (H2). Let us fix a
superlinear function $h:[0,+\infty) \rightarrow \R$ such that
\begin{equation}
\label{suplin}
H(t,q,p) \geq h(|p|_q), \quad \forall (t,q,p) \in [0,1]\times T^*M.
\end{equation}

Let us fix a Riemannian metric on $M$. Such a metric induces metrics on
$TM$ and $T^*M$, together with a horizontal-vertical splitting
\[
TT^*M = T^h T^*M \oplus T^v T^*M,
\]
and isomorphism $T^h_{(q,p)}T^*M \cong T_q M$, $T^v_{(q,p)} T^*M \cong
T_q^*M \cong T_q M$. Let us consider the almost-complex structure $J$
on $T^*M$ which is represented by the matrix
\[
J = \left( \begin{array}{cc} 0 & - I \\ I & 0  \end{array} \right),
\]
with respect to such a splitting. It is $\omega$-compatible, in the
sense that $\omega(J\cdot,\cdot)$ is a Riemannian metric on $T^* M$
(actually, it is the Riemannian metric induced by the one on $M$).

Let us consider the Floer equation
\begin{equation}
\label{floer}
\partial_s u + J(u) [ \partial_t u - X_H(t,u) ] = 0,
\end{equation}
where $u: \R \times [0,1] \rightarrow T^*M$, and $(s,t)$ are the
coordinates on $\R \times [0,1]$. It is a Cauchy-Riemann type first
order elliptic PDE. The solutions of (\ref{floer}) which
do not depend on $s$ are precisely the orbits of the Hamiltonian
vector field $X_H$. If $u$ solves (\ref{floer}), formula
(\ref{fvhaf}) and an integration by parts yield to the energy identity
\[
\int_a^b \int_0^1 |\partial_s u|^2 \, dt \, ds =
\mathbb{A}_H(u(a,\cdot)) - \mathbb{A}_H(u(b,\cdot) + \int_a^b \bigl(
\lambda(u(s,1))[\partial_s u(s,1)] - \lambda(u(s,0))[\partial_s
u(s,0)] \bigr)\, ds.
\]
In particular, if $u$ satisfies the non-local boundary condition
\begin{equation}
\label{hbcu}
(u(s,0),-u(s,1)) \in N^* Q, \quad \forall s\in \R,
\end{equation}
the fact that the Liouville form $\lambda \times \lambda$ of
$T^*(M\times M)$ vanishes on conormal bundles implies that
\begin{equation}
\label{energy}
\int_a^b \int_0^1 |\partial_s u|^2 \, dt \, ds =
\mathbb{A}_H(u(a,\cdot)) - \mathbb{A}_H(u(b,\cdot),
\end{equation}
so the function $s\mapsto \mathbb{A}_H(u(s,\cdot))$ is decreasing.

Let $A\in \R$. We denote by $Z^A(H)$ the space of smooth maps $u:\R
\times [0,1] \rightarrow T^*M$ solving the Floer equation 
(\ref{floer}), satisfying
the boundary condition (\ref{hbcu}), having finite energy,
\[
\int_{-\infty}^{+\infty} \int_0^1  |\partial_s u|^2 \, dt \, ds <
+\infty
\]
and action upper bound
\[
\mathbb{A}_H(u(s,\cdot)) \leq A , \quad \forall s\in \R.
\]
The space $Z^A(H)$, endowed with the $C^{\infty}_{\mathrm{loc}}$ topology,
is a metrizable space.

\begin{lem}
\label{comp}
For every $A\in \R$, the space $Z^A(H)$ is compact.
\end{lem}

The main step in the proof of the above lemma is to show that the
solutions $u\in Z^A(H)$ are uniformly bounded in $C^0$ (see
Theorem 5.4 in \cite{cie94}). Here is where the quadraticity
assumptions (H1') and (H2') are used. Actually, Cieliebak replaces
(H2') by a stronger condition on the second derivatives of $H$. A
proof of the $C^0$ estimate under the above assumptions is
contained in \cite{as06}, Theorem 1.14. Then the facts that
$\omega$ is exact and that the Liouville form $\lambda \times
\lambda$ vanishes on $N^* Q$ allow to exclude bubbling of
$J$-holomorphic spheres or disks, and this yields to uniform $C^1$
bounds for the elements of $Z^A(H)$. The conclusion follows from a
standard elliptic bootstrap. See \cite{cie94}, Theorem 6.3.

We denote by $\overline{H}^*(\cdot,\Z_2)$ the Alexander-Spanier
cohomology with coefficients in $\Z_2$. We denote by $\mathrm{ev}_A$
the evaluation map
\[
\mathrm{ev}_A: Z^A(H) \rightarrow C_Q([0,1],M), \quad  u \mapsto
\pi\circ u(0,\cdot),
\]
where $\pi:T^*M \rightarrow M$ is the projection. The following result
says that if $A$ is large enough then $Z^A(H)$ is non-empty, and that
actually its topology is quite rich.

\begin{lem}
\label{inj}
Let $\alpha$ be a non-zero class in $\overline{H}^*(C_Q([0,1],M),\Z_2)$.
Then there exists $A_0\in \R$,
\[
A_0 = A_0\Bigl( \alpha, h, \max_{(t,q)\in [0,1]\times M} H(t,q,0) \Bigr),
\]
such that for every $A\geq A_0$ the cohomology class $\mathrm{ev}_A^*(\alpha)$
is not zero in $\overline{H}^*(Z^A(H))$.
\end{lem}

This fact is proved in \cite{cie94}, Theorem 7.6. Let us sketch the
argument, in order to show how $A_0$ depends on $H$ and $\alpha$ (this
dependence is not explicitly stated in \cite{cie94}).

Let us fix numbers $S>0$ and $p>1$. Let $\mathscr{B}$ be the Banach
manifold consisting of the maps $u$ in the Sobolev space
$W^{1,p}([-S,S]\times [0,1],T^*M)$ such
that $u(S,t)\in M$ for every $t\in [0,1]$, where $M$ denotes the image
of the zero section in $T^*M$. Let $\mathscr{E}$ be the Banach bundle
over $\mathscr{B}$ whose fiber at $u$ is the space of $L^p$ sections
of $u^*(TT^*M)$. Let $\mathscr{W}$ be the Banach manifold
of paths $\gamma$ in $W^{1-1/p,p}([0,1],M)$ such that
$(\gamma(0),\gamma(1))\in Q$,
and consider the product manifold
$\mathscr{E} \times \mathscr{W}$ as the total space of a bundle over
$\mathscr{B}$. The section
\[
F : \mathscr{B} \rightarrow \mathscr{E} \times \mathscr{W}, \quad u
\mapsto (\partial_s u + J(u) [ \partial_t u - X_H(t,u) ], \pi\circ
u(-S,\cdot)),
\]
is smooth. If $F(u)\in \{0\} \times \mathscr{W}$, then $F$ is a
Fredholm section of index zero at $u$, meaning that the fiberwise
derivative $D^f F(u) : T_u \mathscr{B} \rightarrow \mathscr{E}_u
\times T_{\pi\circ u(-S,\cdot)} \mathscr{W}$
is a Fredholm operator
of index zero. This follows from elliptic estimates and from the fact
that $W^{1-1/p,p}$ is the space of traces of $W^{1,p}$.

Let $K\subset C^{\infty}_Q([0,1],M)$ be a compact set such that
$i^*_K(\alpha)\neq 0$, where $i_K: K\hookrightarrow C_Q([0,1],M)$
denotes the inclusion. If $u=(q,p)\in \mathscr{B}$ is such that
$F(u)\in \{0\} \times K$, then the boundary conditions on $u$ and
(\ref{suplin}) imply the estimates
\begin{eqnarray*}
\mathbb{A}_H(u(-S,\cdot))  =  \int_0^1
p(-S,t)[\partial_t q(-S,t)]\, dt - \int_0^1 H(t,q(-S,t),p(-S,t))\, dt
\\ \leq \int_0^1 \Bigl( \|\partial_t q(-S,\cdot)\|_{\infty} |p(-S,t)|_{q(-S,t)}
 - h\bigl(|p(-S,t)|_{q(-S,t)}\bigr) \Bigr) \,dt, \\
\mathbb{A}_H(u(S,\cdot)) = - \int_0^1 H(t,q(S,t),0)\, dt \geq -  \max_{(t,q)\in
[0,1]\times M} H(t,q,0),
\end{eqnarray*}
Since $q(-S,\cdot)$ belongs to the compact set $K$ and $h$ is
superlinear, we deduce that
\begin{equation}
\label{actbd}
\mathbb{A}_H(u(-S,\cdot)) - \mathbb{A}_H(u(S,\cdot)) \leq A_0,
\end{equation}
where
\begin{equation}
\label{a0}
A_0 := \max_{s\geq 0} \Bigl( \max_{\gamma\in K} \|\gamma'\|_{\infty} s
  - h(s) \Bigr) +  \max_{(t,q)\in [0,1]\times M} H(t,q,0).
\end{equation}
The above action bound and the argument used in the proof of Lemma
\ref{comp} show that $F^{-1}(\{0\}\times K)$ is a compact subset of
$\mathscr{B}$. It can be shown that $Z_{S,K}:=F^{-1}(\{0\}\times K)$ has an
open neighborhood $\mathscr{U}$ in $\mathscr{B}$ such that
\begin{equation}
\label{grado}
\deg ( F|_{\mathscr{U}}, \mathscr{U}, (0,\gamma)) = 1, \quad \forall
\gamma\in K,
\end{equation}
where $\deg$ denotes the Smale $\Z_2$-degree of a Fredholm section of
index zero. Statement (\ref{grado}),
together with the fact that maps between finite dimensional manifolds
of non-zero degree induce injective homomorphisms in cohomology, can
be used to show that $f^*(i_K^*(\alpha))\neq 0$, where $f:Z_{S,K}
\rightarrow K$ is the map $f(u)=\pi\circ u(-S,\cdot)$. If we now let
$S$ tend to infinity, the action bound (\ref{actbd}) implies that the
elements of $Z_{S,K}$ converge to elements of $Z^{A_0}(H)$ in
$C^{\infty}_{\mathrm{loc}}$, and the conclusion of Lemma \ref{inj}
follows by the tautness of the Alexander-Spanier cohomology.

We conclude this section by recalling how Lemmas \ref{comp} and
\ref{inj} imply the existence of solutions of (\ref{he}),
(\ref{hbc}), in the case of Hamiltonians satisfying (H1'), (H2'). Let
$A\in \R$, and consider the continuous flow
\[
(s,u) \mapsto u(s+\cdot,\cdot),
\]
on the compact metrizable space $Z^A(H)$. The continuous function
\[
u \mapsto \mathbb{A}_H(u(0,\cdot))
\]
is a Lyapunov function for such a flow. Then we have the following:

\begin{thm}
\label{lst}
Let $M$ be a compact manifold and let $Q$ be a smooth submanifold of
$M\times M$. Assume that the smooth Hamiltonian $H:[0,1]\times T^*M
\rightarrow \R$ satisfies (H1'), (H2'), and (\ref{suplin}).
\begin{enumerate}
\item Let $\alpha\in \overline{H}^*(C_Q([0,1],M),\Z_2) \cong
  \overline{H}^*(C_Q^{\infty}([0,1],M),\Z_2)$ be a
non-zero Alexander-Spanier cohomology
  class, let $i_K: K \hookrightarrow C^{\infty}_Q([0,1],M)$ be a 
compact inclusion such
  that $i_K^*(\alpha)\neq 0$, and let $A\geq A_0$, where $A_0$ is
  given by (\ref{a0}). Then there is a solution $x$ of (\ref{he}), (\ref{hbc}),
  whose action $\mathbb{A}_H(x)$ equals
\[
c_{\alpha}(H) := \inf_{Z\in \mathscr{Z}_{\alpha}} \max_{u\in Z}
\mathbb{A}_H(u(0,\cdot)),
\]
where $\mathscr{Z}_{\alpha}$ is the set of compact subsets $Z\subset
Z^A(H)$ such that the Alexander-Spanier cohomology class
$\mathrm{ev}_A|_{Z}^*(\alpha)$ is not zero.
\item The number of solutions $x$ of (\ref{he}), (\ref{hbc}), with
  $\mathbb{A}_H(x)\leq A$ is greater or equal than the Alexander-Spanier
  $\Z_2$-cuplength of $Z^A(H)$ plus one.
\item If $\overline{H}^k(C_Q([0,1],M),\Z_2) \neq 0$ for infinitely
many indices
  $k$, then the set of critical levels
\[
\set{c_{\alpha}(H)}{\alpha\in \overline{H}^*(C_Q([0,1],M),\Z_2),
\; \alpha\neq 0}
\]
is not bounded above.
\end{enumerate}
\end{thm}

The first two statements follow from Lusternik-Schnirelman theory
on compact metric spaces, see \cite{hz94}, section 6.3. 
Statement (c) is proved in \cite{cie94}, Proposition 8.4. 

\section{Proof of Theorem \ref{main2}}

Let $H\in C^2([0,1]\times T^*M,\R)$ be a Hamiltonian satisfying
(H1), (H2), (H3). By an easy regularization argument together with
the a priori estimates given by Lemma \ref{hae}, we may assume
that $H$ is smooth. Without loss of generality, we may assume that
the function $a$ appearing in (H1) grows at most linearly, while
the function $h$ appearing in (H2) grows at most as $s\log s$.

\begin{defn}
\label{rm}
We say that $H_0\in C^{\infty}([0,1]\times T^*M,\R)$ is a quadratic
$R$-modification of $H$ if:
\begin{enumerate}
\item $H_0(t,q,p) =  H(t,q,p)$ for $|p|_q \leq R$; 
\item $H_0$ satisfies (H1') and (H2');  
\item $DH_0(t,q,p)[Y(q,p)]-H_0(t,q,p) \geq a(|p|_q)$; 
\item $H_0(t,q,p) \geq h(|p|_q)$.
\end{enumerate}
\end{defn}

It is easy to build a quadratic $R$-modification of a Hamiltonian $H$
satisfying (H1) and (H2), at least for $R$ big enough. Indeed let
us consider a Hamiltonian $H_0$ of the form
\[
H_0(t,q,p):=\varphi(|p|_q)H(t,q,p) + (1-\varphi(|p|_q))C |p|_q^2,
\]
with $C \geq 1$ constant, and $\varphi:[0,+\infty) \to [0,1]$ a
smooth decreasing function taking value $1$ on $[0,R]$
and $0$ on $[R+1,+\infty)$.

Obviously $H_0$ satisfies $(a)$ and $(b)$. Also $(d)$ is satisfied
for $R$ big enough, because $H$ satisfies (H2) and thanks to the
assumption that $h$ grows at most as $s \log s$. Finally, as
$D(|p|_q^2)[Y]=2|p|_q^2$, we have
\[
DH_0[Y]-H_0=\varphi \left(DH[Y]-H \right) + (1-\varphi)C|p|_q^2-
|p|_q \varphi' \left(C|p|_q^2 - H \right).
\]
Since $\varphi'$ is non-positive and it vanishes outside the
interval $[R,R+1]$, it suffices first to take $R$ large enough so
that $|p|_q^2 \geq a(|p|_q)$ for $|p|_q \geq R$ (such $R$ exists
since $c$ grows at most linearly at infinity), and then $C \geq 1$
large enough so that $C|p|_q^2 \geq H(t,q,p)$ for $R \leq
|p|_q\leq R+1$.

\begin{lem}
\label{sqm}
Let $M$ be a compact manifold and $H\in C^{\infty}([0,1]\times
T^*M,\R)$ be a Hamiltonian satisfying (H1), (H2), (H3). For every
$A\in \R$ there exists a number $R(A)$ such that for any $R>R(A)$ and
for any Hamiltonian $H_0$ which is a quadratic $R$-modification of
$H$, the following holds: if $x=(q,p):[0,1]\rightarrow T^*M$ is an
orbit of $X_{H_0}$ such that $\mathbb{A}_{H_0}(x)\leq A$, then
$\|p\|_{\infty} \leq R(A)$.  In particular, such an $x$ is an
orbit of $X_H$, and $\mathbb{A}_{H_0}(x)=\mathbb{A}_{H}(x)$.
\end{lem}

\begin{proof}
Since the Hamiltonian flow of $H$ is globally defined, arguing as
in the proof of Lemma \ref{hae} we find that the set
\[
K :=\set{\phi_t^H \circ (\phi_s^H)^{-1} (q,p)}{(s,t)\in [0,1]^2,
\; (q,p) \in T^*M, \; a(|p|_q) \leq A}
\]
is compact. Let $R(A)$ be the maximum of $|p|_q$ for $(q,p)$ in
$K$. Let $R > R(A)$, $H_0$ a quadratic $R$-modification of $L$,
and $x=(q,p)$ an orbit of $X_{H_0}$ such that
$\mathbb{A}_{H_0}(x)\leq A$. Then there exists $t_0 \in [0,1]$
such that
\[
a(|p(t_0)|) \leq DH_0(t,x(t_0))[Y(x(t_0))]-H_0(t_0,x(t_0)) \leq A,
\]
where the first inequality follows by $(b)$ in Definition
\ref{rm}, which implies that $x(t_0) \in K$, and thus $|p(t_0)|
\leq R(A)$. Let $I \subset [0,1]$ be the maximal interval
containing $t_0$ such that $x(t) \in K$ for every $t \in I$. If $t
\in I$, then $|p(t)| \leq R(A)<R$. Therefore $H$ and $H_0$
coincide in a neighborhood of $(t,x(t))$, and so do the
Hamiltonian vector fields defined by $H$ and $H_0$. Hence, by the
definition of $K$, $x(s)$ belongs to $K$ for $s$ in a neighborhood
of $t$ in $[0,1]$, proving that $I$ is open in $[0,1]$. Being
obviously closed, $I$ coincides with $[0,1]$. Therefore
$\|p\|_{\infty} \leq R(A)$.
\end{proof} \qed

Let $\alpha_1,\dots,\alpha_m$, $m\geq 1$, be cohomology classes of
degree at least one in $H^*(C_Q([0,1],M), \Z_2)$ such that
$\alpha_1 \cup \dots \cup \alpha_m \neq 0$. Since singular
cohomology and Alexander-Spanier cohomology coincide on Banach
manifolds, we may regard $\alpha_1,\dots,\alpha_m$ as
Alexander-Spanier cohomology classes. Let
\[
A = A_0\Bigl(\alpha_1\cup \dots \cup \alpha_m, h, \max_{(t,q)\in
  [0,1]\times M} H(t,q,0) \Bigr),
\]
as given by (\ref{a0}), and let $R>R(A)$, where $R(A)$ is given by
Lemma \ref{sqm}. Let $H_0$ be a quadratic $R$-modification of $H$
(up to taking $R$ larger, we know that a quadratic
$R$-modification always exists). Since $H_0(t,q,p) \geq h(|p|_q)$
and $H_0$ coincides with $H$ on $[0,1]\times M$, Lemma \ref{inj}
implies that
\[
\mathrm{ev}_A^*(\alpha_1) \cup \dots \cup \mathrm{ev}_A^*(\alpha_m) =
\mathrm{ev}_A^*(\alpha_1 \cup \dots \cup \alpha_m) \neq 0 \quad
\mbox{in } \overline{H}^*(Z^A(H_0)).
\]
Therefore, the $\Z_2$-Alexander Spanier cuplength of $Z^A(H_0)$ is at
least $m$. By Theorem \ref{lst} (b), there are at least $m+1$ orbits $x$
of $X_{H_0}$ satisfying the boundary conditions (\ref{hbc}) and the
action estimate $\mathbb{A}_{H_0}(x)\leq A$. By Lemma \ref{sqm} these
curves are also orbits of $X_H$. This proves statement (a) of Theorem
\ref{main2}.

The proof of statement (b) is analogous to the proof of statement
(c) of Theorem \ref{main}, the Morse index being replaced by the
Maslov index. Indeed, the above argument together with the fact
that infinitely many $\Z_2$-homology groups of $C_Q([0,1],M)$ are
non-trivial produces a sequence of solutions of (\ref{he}),
(\ref{hbc}) with diverging Maslov index. See for instance
\cite{sz92} for the definition of the Maslov index and for the
argument yielding to existence of orbits with Maslov index
estimates. By Lemma \ref{hae}, the analogue of Lemma
\ref{boundmorseindex} for the Maslov index holds, so orbits with
diverging Maslov index have diverging action, concluding the proof.

\appendix
\section{Regularity of the Lagrangian action and the Palais-Smale
condition}

Let us say a bit more about the localization argument used in section
\ref{bflt}.
Let $\gamma_0:[0,1]\rightarrow M$ be a continuous curve. Let $U$ be an
open subset of $\R^n$, and let
\begin{equation}
\label{chart2}
[0,1] \times U \rightarrow [0,1] \times M, \quad (t,q) \mapsto
(t,\varphi(t,q)),
\end{equation}
be a smooth coordinate system such that $\gamma_0(t) =
\varphi(t,\tilde\gamma_0(t))$ for every $t\in [0,1]$, for some
continuous curve $\tilde{\gamma}_0:[0,1] \rightarrow U$.
For instance, such a diffeomorphism can be constructed by choosing a
smooth curve $\gamma_1: [0,1]\rightarrow M$ such that $\|\gamma_1 -
\gamma_0\|_{\infty} < \rho$, $\rho$ denoting the injectivity radius of
$M$, and by setting $\varphi(t,q) = \exp_{\gamma_1(t)}[\Psi(t) q]$,
where $q\in \R^n$, $|q|<\rho$, and $\Psi$ is a smooth orthogonal
trivialization of the vector bundle $\gamma_1^*(TM)$ over $[0,1]$.

Let $\Phi$ denote the differential of (\ref{chart2}), that is the
tangent bundle coordinate system
\[
\Phi : [0,1] \times U \times \R^n \rightarrow [0,1] \times TM, \quad
\Phi(t,q,v) = (t,\varphi(t,q),D\varphi(t,q)[(1,v)]).
\]
Up to replacing $U$ by a smaller open set, we may assume that $\Phi$
is bounded together with its inverse. Then the map
\[
\varphi_* : W^{1,2}([0,1],U) \rightarrow W^{1,2}([0,1],M),  \quad
\varphi_* (\tilde{\gamma})(t) = \varphi(t,\gamma(t)),
\]
is a smooth local coordinate system on the Hilbert manifold
$W^{1,2}([0,1],M)$, such that $D\varphi_*$ is bounded together with its
inverse from the standard metric of $W^{1,2}([0,1],\R^n)$ to the Riemannian
metric on $W^{1,2}([0,1],M)$ defined in (\ref{metric}). If moreover
the curve $\gamma_0:[0,1]\rightarrow M$ is in $W^{1,2}$, then so is
$\tilde{\gamma}_0$, and $\varphi_*(\tilde{\gamma}_0) = \gamma_0$.

The Lagrangian $L$ on $[0,1]\times TM$ can be pulled back onto
$[0,1]\times U \times \R^n$ by setting $\tilde{L} = L \circ
\Phi$. Notice that conditions (L1') and (L2') are invariant with
respect to coordinates transformations of the form $\Phi$ (up to
changing the constants), so $\tilde{L}$ satisfies these
conditions if $L$ does. Since
$L(t,\varphi_*(\tilde\gamma)(t),\varphi_*(\tilde\gamma)'(t)) =
\tilde{L}(t,\tilde\gamma(t),\tilde\gamma'(t))$, we have
\[
\mathbb{A}_L(\varphi_*(\tilde{\gamma})) =
\mathbb{A}_{\tilde{L}}(\tilde\gamma).
\]
These facts allow to reduce all the arguments on $W^{1,2}([0,1],M)$
which are $C^0$-local, in the sense that they involve only a
$C^0$-neighborhood of some curve, to the case of a Lagrangian on an
open subset of the Euclidean space.

\paragraph{Proof of Proposition \ref{reg}.} By the above localization
arguments, we may assume that $L$ is a Lagrangian on $[0,1]\times U
\times \R^n$, with $U$ an open subset of $\R^n$.
If $\gamma\in W^{1,2}([0,1],U)$, $\xi\in W^{1,2}([0,1],\R^n)$, and
$h\in \R$ has a small absolute value, we have
\begin{equation}
\label{regC1}
\frac{1}{h} \Bigl(\mathbb A_L(\gamma+h \xi) -\mathbb
  A_L(\gamma)\Bigr) = \int_0^1 \int_0^1 \Bigl( D_q
  L(t,\gamma+ h s \xi, \gamma'+ h s \xi') [\xi]
+ D_v L(t,\gamma+ h s \xi,\gamma'+h s \xi') [\xi'] \Bigr) \,dt \, ds.
\end{equation}
Assumption (L2') implies that
\begin{equation}
\label{bfo}
|D_v L(t,q,v)| \leq \ell_2 (1+|v|_q), \quad |D_q L(t,q,v)| \leq \ell_2
(1+|v|_q^2),
\end{equation}
for some constant $\ell_2$. These bounds and the dominated convergence
theorem show that the quantity (\ref{regC1}) converges to
\[
d\mathbb{A}_L(\gamma)[\xi] = \int_0^1 \bigl( D_q
  L(t,\gamma,\gamma')[\xi] - D_v L(t,\gamma,\gamma')[\xi'] \bigr)\,
dt
\]
for $h\rightarrow 0$. Since $d\mathbb{A}_L(\gamma)$ is a bounded
linear functional on $W^{1,2}([0,1],\R^n)$, $\mathbb A_L$ is Gateaux
differentiable and $d\mathbb{A}_L(\gamma)$ is its Gateaux
differential at $\gamma$. Since
$d\mathbb{A}_L(\gamma)$ continuously depends on $\gamma\in
W^{1,2}([0,1],U)$, the total differential theorem implies that
$\mathbb{A}_L$ is continuously Fr\'ech\'et differentiable, with
$D\mathbb{A}_L = d\mathbb{A}_L$.

Similarly, by (L2') and by the dominated convergence theorem, the quantity
\begin{equation}
\label{regC2}
\begin{split}
\frac{1}{h} \bigl(&D\mathbb A_L(\gamma+h \eta)[\xi] -
D\mathbb A_L(\gamma)[\xi]\bigr)\\ =&
\int_0^1 \int_0^1 \Bigl( D_{vv} L(t,
\gamma+ h s \eta, \gamma'+ h s \eta')[\xi',\eta'] +
D_{vq} L(t,\gamma+ h s \eta, \gamma'+h s \eta')[\xi',\eta]\\
&+ D_{qv} L(t,\gamma+h s \eta, \gamma'+h s \eta')[\xi,\eta']
+ D_{qq} L(t,\gamma+h s \eta, \gamma'+h s \eta')[\xi,\eta]\Bigr) \, dt\,ds,
\end{split}
\end{equation}
converges to
\begin{eqnarray*}
d^2 \mathbb{A}_L (\gamma)[\xi,\eta] =
  \int_0^1 \Bigl( D_{vv} L(t,\gamma,\gamma')
[\xi',\eta'] + D_{qv} L(t,\gamma,\gamma') [\xi,\eta'] \\
+ D_{vq}
  L(t,\gamma,\gamma') [\xi',\eta] + D_{qq} L(t,\gamma,\gamma') [\xi,\eta]
  \Bigr) \, dt,
\end{eqnarray*}
for $h\rightarrow 0$. Since $d^2 \mathbb{A}_L(\gamma)$ is a bounded
symmetric bilinear form on $W^{1,2}([0,1],\R^n)$ which
depends continuously on
$\gamma$, we conclude that $\mathbb{A}_L$ is $C^2$, and that its
second differential is $D^2 \mathbb{A}_L = d^2
\mathbb{A}_L(\gamma)$. Since the second differential at critical
points is invariantly defined, in the sense that
\[
D^2
\mathbb{A}_L(\psi(\gamma))[D\psi(\gamma)\xi,D\psi
(\gamma)\eta] = D^2 (\mathbb{A}_{L}\circ \psi) (\gamma) [\xi,\eta],
\]
this concludes the proof of Proposition
\ref{reg}.

\paragraph{Proof of Proposition \ref{ps}.}
Let $(\gamma_h)$ be a sequence in $W_Q^{1,2}([0,1],M)$ such that
$\mathbb A_L(\gamma_h)$ is bounded and $\|D\mathbb
A_L(\gamma_h)|_{T_{\gamma} W^{1,2}_Q([0,1],M)}\|_*$
is infinitesimal. We have to show that
$(\gamma_h)$ is compact in $W_Q^{1,2}([0,1],M)$.
First we observe that by assumption (L1'), $L(t,q,v) \geq \ell_0
|v|_q^2/2 - C$, for some constant $C$. Together with the bound on
the action, this implies that the sequence $(\gamma_h')$ is bounded
in $L^2$. Therefore,
\[
\dist (\gamma_h(t),\gamma_h(s)) \leq \int_s^t |\gamma_h'(\sigma)|\, d
\sigma \leq |s-t|^{1/2} \left( \int_0^1 |\gamma_h'(\sigma)|^2\, d
\sigma \right)^{1/2},
\]
so $(\gamma_h)$ is equi-1/2-H\"older continuous,
and up to a subsequence we may assume that it converges uniformly
to some $\gamma\in C([0,1],M)$ with $(\gamma(0),\gamma(1))\in Q$.
If $\varphi_*$ is a smooth coordinate
system on $W^{1,2}([0,1],M)$ constructed as above with
$\gamma_0:=\gamma$, $\gamma_h$ eventually belongs to the image of
$\varphi_*$, so we may assume that $L$ is a Lagrangian on $[0,1]\times U
\times \R^n$ and $\gamma_h\in W^{1,2}([0,1],U)$, where
$U$ is an open subset of $\R^n$. Up to choosing the diffeomorphism
$(t,q) \mapsto (t,\varphi(t,q))$ properly, we may also assume that there is
an affine subspace $V=\zeta + V_0$ of $\R^n \times \R^n$ such that
\[
[\varphi(0,\cdot) \times \varphi(1,\cdot)]^{-1} Q = (U \times U) \cap
V.
\]
Notice also that, since the coordinate system
$\varphi_*$ is $C^1$-bounded together with its inverse,
$D\mathbb{A}_L(\gamma_h)$ converges to zero strongly in the dual of
\[
W^{1,2}_{V_0}([0,1],\R^n) := \set{\xi\in
  W^{1,2}([0,1],\R^n)}{(\xi(0), \xi(1)) \in V_0}.
\]
The fact that $(\gamma_h')$ is bounded
in $L^2$ now implies that $\gamma\in W^{1,2}([0,1],U)$, and up to a
subsequence, $(\gamma_h)$ converges weakly to $\gamma$ in
$W^{1,2}([0,1],\R^n)$ and strongly in $L^2$.
We must show that this convergence is strong in $W^{1,2}$.
Since $D\mathbb{A}_L(\gamma_h)$ is infinitesimal in the dual
of $W^{1,2}_{V_0}([0,1],\R^n)$,
$(\gamma_h)$ is bounded in $W^{1,2}$, and $\gamma_h-\gamma\in
W^{1,2}_{V_0}([0,1],\R^n)$,
the sequence $D\mathbb A_L(\gamma_h)[\gamma_h-\gamma]$ is
infinitesimal, that is
\[
\int_0^1  D_q L(t,\gamma_h,\gamma_h')
[\gamma_h-\gamma]\, dt  + \int_0^1 D_v
L(t,\gamma_h,\gamma_h')[\gamma_h'-\gamma']
\, dt \rightarrow 0 \quad \mbox{for } h \rightarrow \infty.
\]
Since $D_q L(t,\gamma_h,\gamma_h')$ is bounded
in $L^2$ by (\ref{bfo}), and $\gamma_h-\gamma$ converges strongly to $0$
in $L^2$, the first integral in the above expression tends to zero, so
\begin{equation}
\label{inf}
\int_0^1 D_v L(t,\gamma_h,\gamma_h')[\gamma_h'-\gamma']
\, dt \rightarrow 0 \quad \mbox{for } h \rightarrow \infty.
\end{equation}
For a.e.\ $t\in [0,1]$ we have by (L1')
\begin{eqnarray*}
D_v L(t,\gamma_h,\gamma_h')[\gamma_h' - \gamma_h] -D_v
L(t,\gamma_h,\gamma')[\gamma_h' - \gamma_h] \\ = \int_0^1 D_{vv}
L(t,\gamma_h,\gamma' + s(\gamma_h' - \gamma')) [\gamma_h' -
\gamma',\gamma_h' - \gamma'] \, ds   \geq \ell_0 |\gamma_h'(t) -
\gamma'(t)|^2.
\end{eqnarray*}
Integrating this inequality over $[0,1]$ we get
\begin{equation}
\label{dis}
\ell_0 \int_0^1 |\gamma_h' - \gamma'|^2 \, dt \leq \int_0^1 D_v
L(t,\gamma_h,\gamma_h') [\gamma_h'-\gamma'] \, dt - \int_0^1 D_v
L(t,\gamma_h,\gamma') [\gamma_h'-\gamma']\, dt.
\end{equation}
The first integral on the right-hand side tends to $0$ for
$h\rightarrow \infty$ by (\ref{inf}). Since $D_v
L(\cdot,\gamma_h,\gamma')$ converges strongly to $D_v
L(\cdot,\gamma,\gamma')$ in $L^2$ by (\ref{bfo}), and since
$(\gamma_h'-\gamma')$ converges weakly to $0$ in $L^2$, the last
integral in (\ref{dis}) is also infinitesimal. Therefore,
(\ref{dis}) implies that $(\gamma_h)$ converges strongly to
$\gamma$ in $W^{1,2}([0,1],\R^n)$, concluding the
proof.

\providecommand{\bysame}{\leavevmode\hbox to3em{\hrulefill}\thinspace}
\providecommand{\MR}{\relax\ifhmode\unskip\space\fi MR }
\providecommand{\MRhref}[2]{%
  \href{http://www.ams.org/mathscinet-getitem?mr=#1}{#2}
}
\providecommand{\href}[2]{#2}

\end{document}